\documentclass[a4paper,9pt]{article}
\usepackage{amsmath, amsfonts, amssymb}
\usepackage[english, russian]{babel}
\usepackage[cp1251]{inputenc}

\usepackage{graphicx}
\numberwithin{equation}{section}
\begin{document}
\sloppy
\begin{center}
\textbf{Fundamental solutions for a class of multidimensional  elliptic \\ equations with several singular coefficients}\\[5pt]
% Авторы
\textbf{Ergashev T.G.\\}
% обязательное поле!
\medskip
{ Institute of Mathematics, Uzbek
Academy of Sciences,  Tashkent, Uzbekistan. \\

{\verb ergashev.tukhtasin@gmail.com }}\\
\end{center}

\begin{quote} % Аннотация на английском

The main result of the present paper is the construction of
fundamental solutions for a class of multidimensional elliptic
equations with several singular coefficients. These fundamental
solutions  are directly connected with  multiple hypergeometric
functions and  the decomposition formula is required for their
investigation which would express the multivariable hypergeometric
function in terms of products of several simpler hypergeometric
functions involving fewer variables. In this paper, such a formula
is proved instead of a previously existing recurrence formula.The
order of singularity and other properties of the fundamental
solutions that are necessary for solving boundary
value problems for degenerate second-order elliptic equations are determined.
\end{quote}

\section{Introduction}

It is known that fundamental solutions have an essential role in
studying partial differential equations. Formulation and solving
of many local and non-local boundary value problems are based on
these solutions. Moreover,\,\,\,fundamental solutions appear as
potentials, for instance, as simple-layer and double-layer
potentials in the theory of potentials.

The explicit form of fundamental solutions gives a possibility to
study the considered equation in detail. For example, in the works
of Barros-Neto and Gelfand \cite{{BN1},{BN2},{BN3}} fundamental
solutions for Tricomi operator, relative to an arbitrary point in
the plane were explicitly calculated. We also mention Leray's work
\cite{L}, which it was described as a general method, based upon
the theory of analytic functions of several complex variables, for
finding fundamental solutions for a class of hyperbolic linear
differential operators with analytic coefficients. In this
direction we would like to note the works  \cite{{HK},{I}}, where
three-dimensional fundamental solutions for elliptic equations
were found. In the works \cite{{EH},{M},{U}}, fundamental
solutions for a class of multidimensional degenerate elliptic
equations with spectral parameter were constructed. The found
solutions can be applied to solving some boundary value problems
\cite{{Erg},{GC},{Kar},{KN},{SH7},{SH8},{SHCh}}.

Let us consider the generalized Helmholtz equation with a several
singular coefficients
\begin{equation}\label{e1}
L_{(\alpha)}^{m}  (u):={\sum\limits_{i = 1}^{m} {u_{x_{i} x_{i}}}
} + {\sum\limits_{j = 1}^{n} {{\frac{{2\alpha _{j}}} {{x_{j}}} }}}
u_{x_{j}} = 0
\end{equation}
in the domain $R^{n+}_m:=\left\{
\left(x_1,...,x_m\right):x_1>0,...,x_n>0\right\}\,$, where $m$ is
a dimension of the Euchlidean space, $n$ is a number of the
singular coefficients of equation (\ref{e1}); $m \ge 2,0 \le n \le
m$; $\alpha_j$ are real constants and $0<2\alpha_j<1, $
$j=1,...,n$; $(\alpha)=(\alpha_1,...,\alpha_n) $.

Various modifications of the equation (\ref{e1}) in the two- and
three-dimensional cases were considered in many papers
\cite{{A},{F},{H},{HK},{K},{MC},{R},{UrK},{Wt},{Wn4},{Wn5}}.

Fundamental solutions for elliptic equations with singular
coefficients are directly connected with hypergeometric functions.
Therefore, basic properties such as decomposition formulas,
integral representations, formulas of analytical continuation,
formulas of differentiation for hypergeometric functions are
necessary for studying fundamental solutions.

Since the aforementioned properties of hypergeometric functions of
Gauss, Appell, Kummer were known \cite{AP}, results on
investigations of elliptic equations with one or two singular
coefficients were successful. In the paper \cite{HK} when finding
and studying the fundamental solutions of equation (\ref{e1}) for
$m = 3$, an important role was played the decomposition formula of
Hasanov and Srivastava \cite{{HS6},{HS7}}, however, the recurrence
of this formula did not allow further advancement in the direction
of increasing the number of singular coefficients.

In the present paper we construct the $2^n$ fundamental solutions
for equation (\ref{e1}) in an explicit form and we prove a new
formula for the expansion of several Lauricella hypergeometric
functions by simple Gauss, with which it is possible to reveal
that the found hypergeometric functions have a singularity of
order $1/r^{m-2}$ at $r \to 0.$

The plan of this paper is as follows. In Section 2 we briefly give
some preliminary information, which will be used later.We
transform the recurrence decomposition formula of Hasanov and
Srivastava \cite{HS6} to the form convenient for further research.
Also some constructive formulas for the operator   $L$ are given.
In Section 3 we describe the method of finding fundamental
solutions for the considered equation and in section we show what
order of singularity the found solutions will have.

 \section{Preliminaries}

Below we give definition of  Pochhammer symbol and some formulas
for Gauss hypergeometric functions of one and two variables,
Lauricella hypergeometric functions of three and more variables,
which will be used in the next section.

A symbol $\left( {\kappa}  \right)_{\nu}  $ denotes the general
Pochhammer symbol or the shified factorial, since $\left( {1}
\right)_{l} = l!$ $ \left( {l  \in N \cup \{0\};\,\,N: =
\{1,2,3,...\}} \right),$ which is defined (for $\kappa ,\nu \in
C)$, in terms of the familiar Gamma function, by
$$
(\kappa )_{\nu}  : = {\frac{{\Gamma (\kappa + \nu )}}{{\Gamma
(\kappa )}}} = {\left\{ {{\begin{array}{*{20}c}

{1\,\,\,\,\,\,\,\,\,\,\,\,\,\,\,\,\,\,\,\,\,\,\,\,\,\,\,\,\,\,\,\,\,\,\,\,\,\,\,\,\,\,\,\,\,(\nu
= 0;\kappa \in C\backslash \{0\})} \hfill \\
 {\kappa (\kappa + 1)...(\kappa + l - 1)\,\,\,\,\,(\nu = l \in N;\kappa \in
C),} \hfill \\
\end{array}}}  \right.}
$$

\noindent it being understood conventionally that $\left( {0}
\right)_{0} : = 1$ and assumed tacitly that the $\Gamma - $
quotient exists.

A function
$$F\left[
{{\begin{array}{*{20}c}
 {a,b;} \hfill \\
 {c;} \hfill \\
\end{array}} x}  \right]=\sum\limits_{k=0}^\infty\frac{(a)_k(b)_k}{k!(c)_k}x^k,\,\,|x|<1$$
is known as the Gauss hypergeometric function and an equality
\begin{equation}\label{e21}
F\left[ {{\begin{array}{*{20}c}
 {a,b;} \hfill \\
 {c;} \hfill \\
\end{array}} 1}  \right]=\frac{\Gamma(c)\Gamma(c-a-b)}{\Gamma(c-a)\Gamma(c-b)}, c\neq0,-1,-2,..., \textrm{Re}(c-a-b)>0
\end{equation} holds \cite[ p.73,(14)]{Erd}. Moreover, the following autotransformer
formula \cite[ p.76,(22)]{Erd}
\begin{equation}\label{e22}F\left[
{{\begin{array}{*{20}c}
 {a,b;} \hfill \\
 {c;} \hfill \\
\end{array}} x}  \right]=(1-x)^{-b}F\left[
{{\begin{array}{*{20}c}
 {c-a,b;} \hfill \\
 {c;} \hfill \\
\end{array}} \frac{x}{x-1}}  \right] \end{equation}
is valid.

 The hypergeometric function of
$n$ variables has a form \cite[pp.114-115(1),(5)]{AP} (see also
\cite[p.33]{SK})
\begin{equation}\label{e23}F^{(n)}_A \left[ {{\begin{array}{*{20}c}
 {a,b_1,...,b_n;} \hfill \\
 {c_1,...,c_n;} \hfill \\
\end{array}} x_1,...,x_n}  \right]
=\sum\limits_{m_1,...,m_n=0}^\infty\frac{(a)_{m_1+...+m_n}(b_1)_{m_1}...(b_n)_{m_n}}{m_1!...m_n!(c_1)_{m_1}...(c_n)_{m_n}}x_1^{m_1}...x_n^{m_n},
\end{equation}
where  $|x_1|+...+|x_n|<1,  n \in {\rm N}$.

For a given multivariable function, it is useful to fund a
decomposition formula which would express the multivariable
hypergeometric function in terms of products of several simpler
hypergeometric functions involving fewer variables.

In the case of two variables  for the function
\begin{equation*}F_2 \left[ {{\begin{array}{*{20}c}
 {a,b_1,b_2;} \hfill \\
 {c_1,c_2;} \hfill \\
\end{array}} x,y}  \right]
=\sum\limits_{i,j=0}^\infty\frac{(a)_{i+j}(b_1)_{i}(b_2)_{j}}{i!j!(c_1)_{i}(c_2)_{j}}x^{i}y^{j}
\end{equation*}
was known expansion formula \cite{BC1}
\begin{equation}\label{e25}F_2 \left[ {{\begin{array}{*{20}c}
 {a,b_1,b_2;} \hfill \\
 {c_1,c_2;} \hfill \\
\end{array}} x,y}  \right]
=\sum\limits_{k=0}^\infty\frac{(a)_k(b_1)_k(b_2)_k}{k!(c_1)_{k}(c_2)_{k}}x^{k}y^{k}F\left[
{{\begin{array}{*{20}c}
 {a+k,b_1+k;} \hfill \\
 {c_1+k;} \hfill \\
\end{array}} x}  \right]F\left[ {{\begin{array}{*{20}c}
 {a+k,b_2+k;} \hfill \\
 {c_2+k;} \hfill \\
\end{array}} y}  \right].
\end{equation}

Following the works  \cite{{BC1},{BC2}} Hasanov and Srivastava
\cite{HS6} found following decomposition formula for the
Lauricella function of three variables
\begin{equation}\label{e26}\begin{array}{l}F^{(3)}_A \left[ {{\begin{array}{*{20}c}
 {a,b_1,b_2,b_3;} \hfill \\
 {c_1,c_2,c_3;} \hfill \\
\end{array}} x,y,z}  \right]
={\sum\limits_{i,j,k=0}^\infty\frac{(a)_{i+j+k}(b_1)_{j+k}(b_2)_{i+k}(b_3)_{i+j}}{i!j!k!(c_1)_{j+k}(c_2)_{i+k}(c_3)_{i+j}}}
\\\cdot
x^{j+k}y^{i+k}z^{i+j}
  F\left[ {{\begin{array}{*{20}c}
 {a+j+k,b_1+j+k;} \hfill \\
 {c_1+j+k;} \hfill \\
\end{array}} x}  \right]\\
\cdot F\left[ {{\begin{array}{*{20}c}
 {a+i+j+k,b_2+i+k;} \hfill \\
 {c_2+i+k;} \hfill \\
\end{array}} y}  \right]F\left[ {{\begin{array}{*{20}c}
 {a+i+j+k,b_3+i+j;} \hfill \\
 {c_3+i+j;} \hfill \\
\end{array}} z}  \right]\end{array}
\end{equation}
and they proved that for all $n\in N\backslash\{1\}$ is true the
recurrence formula \cite {HS6}
\begin{equation}
\label{e27} \begin{array}{l} F_{A}^{(n)}\left[
{{\begin{array}{*{20}c}
 {a,b_1,...,b_n;} \hfill \\
 {c_1,...,c_n;} \hfill \\
\end{array}} x_1,...,x_n}  \right]
 = {\sum\limits_{m_{2} ,...,m_{n} = 0}^{\infty}  {{\frac{{(a)_{m_{2} + \cdot
\cdot \cdot + m_{n}}  (b_{1} )_{m_{2} + \cdot \cdot \cdot + m_{n}}
(b_{2} )_{m_{2}}  \cdot \cdot \cdot (b_{n} )_{m_{n}}} } {{m_{2} !
\cdot \cdot \cdot m_{n} !(c_{1} )_{m_{2} + \cdot \cdot \cdot +
m_{n}}  (c_{2} )_{m_{2}}  \cdot \cdot \cdot (c_{n} )_{m_{n}}} }
}}}
\\ \cdot x_{1}^{m_{2} + \cdot \cdot \cdot + m_{n}} x_{2}^{m_{2}}
\cdot \cdot \cdot x_{n}^{m_{n}}    F\left[ {{\begin{array}{*{20}c}
 {a + m_{2} + \cdot \cdot \cdot + m_{n},b_{1} + m_{2} +
\cdot \cdot \cdot
+ m_{n};} \hfill \\
 {c_{1} + m_{2} +
\cdot \cdot \cdot
+ m_{n};} \hfill \\
\end{array}} x_1}  \right] \\
 \cdot F_{A}^{(n - 1)} \left[ {{\begin{array}{*{20}c} {a + m_{2} + \cdot \cdot \cdot + m_{n} ,b_{2}
+ m_{2} ,...,b_{n} + m_{n} ;}\hfill \\
c_{2} + m_{2} ,....,c_{n} + m_{n}
;\hfill \\ \end{array}} x_{2} ,...,x_{n}}  \right]. \\
\end{array}
\end{equation}

Further study of the properties of the Lauricella function
(\ref{e23}) showed that the formula (\ref{e27}) can be reduced to
a more convenient form.

{\textbf{Lemma.} \textit{The following formula holds true at $n\in
N$
\begin{equation}
\label{e28}\begin{array}{l} {F_{A}^{(n)} {\left[
{{\begin{array}{*{20}c}{a,b_{1}
,....,b_{n} ;} \hfill \\ {c_{1} ,....,c_{n} ;} \hfill \\
\end{array}} x_{1} ,...,x_{n}}  \right]}
 = {\sum\limits_{{\mathop {m_{i,j} = 0}\limits_{(2 \le i \le j \le n)}
}}^{\infty}  {{\frac{{(a)_{N_{2} (n,n)}}} {{{\mathop {m_{2,2}
!m_{2,3} ! \cdot \cdot \cdot m_{i,j} ! \cdot \cdot \cdot m_{n,n}
!}\limits_{(2 \le i \le j \le n)}}} }}}}} \hfill\\
 \cdot{\prod\limits_{k = 1}^{n} {{ {{\frac{{(b_{k} )_{M_{2} (k,n)}
}}{{(c_{k} )_{M_{2} (k,n)}}} x_{k}^{M_{2} (k,n)}
F\left[{\begin{array}{*{20}c} {a + N_{2} (k,n),b_{k} + M_{2}
(k,n);} \hfill \\ c_{k} + M_{2} (k,n);\hfill \\ \end{array} x_{k}}
\right]} } }}},
\end{array}
\end{equation}
where
$$
M_{l} (k,n) = {\sum\limits_{i = l}^{k} {m_{i,k} +}}
{\sum\limits_{i = k + 1}^{n} {m_{k + 1,i}}},\,\, \quad N_{l} (k,n)
= {\sum\limits_{i = l}^{k + 1} {{\sum\limits_{j = i}^{n}
{m_{i,j}}} } } , \, l \in N.
$$
}

\textbf{Proof}. We carry out the proof by the method mathematical
induction.

The equality (\ref{e28}) in the case $n=1$ is obvious.

Let $n = 2$. Since $M_2(1,2)=M_2(2,2)=N_2(1,2)=N_2(2,2)=m_{2,2},$
we obtain the formula (\ref{e25}).

For the sake of interest, we will check the formula (\ref{e28}) in
yet another value of $n$.

Let $n=3.$ In this case
$$M_2(1,3)=m_{2,2}+m_{2,3},\,\, M_2(2,3)=m_{2,2}+m_{3,3},\,\, M_2(3,3)=m_{2,3}+m_{3,3},$$
$$N_2(1,3)=m_{2,2}+m_{2,3},\,\, N_2(2,3)= N_2(3,3)=m_{2,2}+m_{2,3}+m_{3,3}.$$
For brevity, making the substitutions
$m_{2,2}:=i,\,\,m_{2,3}:=j,\,\,m_{3,3}:=k$, we obtain the formula
(\ref{e26}).

So the formula (\ref{e28}) works for $n=1,$ $n=2$ and $n=3$.

Now we assume that for $n = s$ equality (\ref{e28}) holds; that
is, that
\begin{equation}
\label{e29}
\begin{array}{l}
 F_{A}^{(s)} {\left[ \begin{array}{*{20}c}{a,b_{1} ,....,b_{s} ;}\hfill\\ c_{1} ,....,c_{s} ; \hfill\\ \end{array} x_{1}
,...,x_{s}  \right]}
 = {\sum\limits_{{\mathop {m_{i,j} = 0}\limits_{(2 \le i \le j \le s)}
}}^{\infty}  {{\frac{{(a)_{N_{2} (s,s)}}} {{{\mathop {m_{2,2}
!m_{2,3} ! \cdot \cdot \cdot m_{i,j} ! \cdot \cdot \cdot m_{s,s}
!}\limits_{(2 \le i
\le j \le s)}}} }}}}  \\
 \cdot {\prod\limits_{k = 1}^{s} {{{{\frac{{(b_{k} )_{M_{2} (k,s)}
}}{{(c_{k} )_{M_{2} (k,s)}}} }x_{k}^{M_{2} (k,s)} F\left[
\begin{array}{*{20}c}{a + N_{2} (k,s),b_{k} + M_{2} (k,s);}\hfill\\ c_{k} +
M_{2} (k,s);\hfill\\ \end{array}  x_{k} \right]} }}} .
\\
 \end{array}
\end{equation}

Let $n=s+1.$ We prove that  following formula
\begin{equation}
\label{e210}
\begin{array}{l}
 F_{A}^{(s + 1)} {\left[ \begin{array}{*{20}c}{a,b_{1} ,....,b_{s+1} ;}\hfill\\ c_{1} ,....,c_{s+1} ; \hfill\\ \end{array} x_{1}
,...,x_{s+1}  \right]}
 = {\sum\limits_{{\mathop {m_{i,j} = 0}\limits_{(2 \le i \le j \le s + 1)}
}}^{\infty}  {{\frac{{(a)_{N_{2} (s + 1,s + 1)}}} {{{\mathop
{m_{2,2} !m_{2,3} ! \cdot \cdot \cdot m_{i,j} ! \cdot \cdot \cdot
m_{s + 1,s + 1}
!}\limits_{(2 \le i \le j \le s + 1)}}} }}}}  \\
 \cdot {\prod\limits_{k = 1}^{s + 1} {{ {{\frac{{(b_{k} )_{M_{2} (k,s + 1)}
}}{{(c_{k} )_{M_{2} (k,s + 1)}}} }x_{k}^{M_{2} (k,s + 1)} F\left[
{{\begin{array}{*{20}c}
 {a + N_{2} (k,s + 1),b_{k} + M_{2} (k,s + 1);} \hfill \\
 {c_{k} + M_{2} (k,s + 1);} \hfill \\
\end{array}} x_{k}}  \right]} }}}  \\
 \end{array}
\end{equation}
is valid.

We write the Hasanov-Srivastava's formula (\ref{e27}) in the form

\begin{equation}
\label{e211}
\begin{array}{l}
 F_{A}^{(s + 1)} {\left[ {{\begin{array}{*{20}c}{a,b_{1} ,....,b_{s + 1} ;}\hfill \\ { c_{1} ,....,c_{s +
 1};}\hfill \\
 \end{array}}
x_{1} ,...,x_{s + 1}}  \right]} \\
 = {\sum\limits_{m_{2,2} ,...,m_{2,s + 1} = 0}^{\infty}  {{\frac{{(a)_{N_2(1,s+1)}  (b_{1} )_{M_2(1,s+1)}  (b_{2} )_{m_{2,2}}  \cdot \cdot \cdot (b_{s + 1}
)_{m_{2,s + 1}}} } {{m_{2,2} ! \cdot \cdot \cdot m_{2,s + 1}
!(c_{1} )_{M_2(1,s+1)}  (c_{2} )_{m_{2,2}} \cdot \cdot \cdot (c_{s
+ 1})_{m_{2,s + 1}}} } }}}   \\
 \cdot x_{1}^{M_2(1,s+1)}  x_{2}^{m_{2,2}} \cdot
\cdot \cdot x_{s + 1}^{m_{2,s + 1}} F\left[\begin{array}{*{20}c}{a + N_2(1,s+1) ,b_{1} + M_2(1,s+1) ;} \hfill\\ c_{1} +M_2(1,s+1);\hfill\\ \end{array} x_{1} \right] \\
 \cdot F_{A}^{(s)} {\left[ {{\begin{array}{*{20}c}
 {a + N_2(1,s+1) ,b_{2} + m_{2,2} ,...,b_{s + 1} +
m_{2,s + 1} ;} \hfill \\
 {c_{2} + m_{2,2} ,....,c_{s + 1} + m_{2,s + 1}};  \hfill \\
\end{array}} x_{2} ,...,x_{s + 1}}  \right]}. \\
 \end{array}
\end{equation}

By virtue of  the formula (\ref{e29}) we have

\begin{equation}
\label{e212}
\begin{array}{l}
 F_{A}^{(s)} \left[ \begin{array}{*{20}c}{a + N_2(1,s+1) ,b_{2} + m_{2,2} ,...,b_{s +
1} + m_{2,s + 1} ;}\hfill\\ c_{2} + m_{2,2} ,...,c_{s + 1} +
m_{2,s + 1} ;\hfill\\ \end{array} x_{2} ,...,x_{s +
1}  \right] \\
 = {\sum\limits_{{\mathop {m_{i,j} = 0}\limits_{(3 \le i \le j \le s + 1)}
}}^{\infty}  {{\frac{{\left( {a +N_2(1,s+1)} \right)_{N_{3} (s +
1,s + 1)}}} {{{\mathop {m_{3,3} !m_{3,4} ! \cdot \cdot \cdot
m_{i,j} ! \cdot \cdot \cdot m_{s + 1,s + 1} !}\limits_{(3 \le i
\le j \le s + 1)} }}}}}} \prod\limits_{k = 2}^{s + 1}
{\frac{{(b_{k} + m_{2,k} )_{M_{3} (k,s + 1)}}} {{(c_{k} + m_{2,k}
)_{M_{3} (k,s + 1)}}}
}x_{k}^{M_{3} (k,s + 1)}  \\
 \cdot F\left[ {{\begin{array}{*{20}c}
 {a + N_2(1,s+1) + N_{3} (k,s + 1),b_{k} + m_{2,k} + M_{3} (k,s +
1);} \hfill \\
 {c_{k} + m_{2,k} + M_{3} (k,s + 1);} \hfill \\
\end{array}} x_{k}}  \right]  . \\
 \end{array}
\end{equation}

Substituting from (\ref{e212}) into (\ref{e211}) we obtain
\begin{equation*}
\begin{array}{l}
 F_{A}^{(s + 1)} {\left[ {a,b_{1} ,....,b_{s + 1} ;c_{1} ,....,c_{s + 1}
;x_{1} ,...,x_{s + 1}}  \right]} \\
 = {\sum\limits_{{\mathop {m_{i,j} = 0}\limits_{(2 \le i \le j \le s + 1)}
}}^{\infty}  {{\frac{{\left( {a} \right)_{N_2(1,s+1) + N_{3} (s +
1,s + 1)}}} {{{\mathop {m_{2,2} !m_{2,3} ! \cdot \cdot \cdot
m_{i,j} ! \cdot \cdot \cdot m_{s + 1,s + 1} !}\limits_{(2 \le i
\le j \le s + 1)}}} }}}} \prod\limits_{k = 1}^{s + 1}
{\frac{{(b_{k} )_{m_{2,k} + M_{3} (k,s + 1)}}} {{(c_{k} )_{m_{2,k}
+ M_{3} (k,s + 1)}}}
}x_{k}^{m_{2,k} + M_{3} (k,s + 1)}   \\
\cdot F\left[ {{\begin{array}{*{20}c} {a + N_2(1,s+1) + N_{3} (k,s
+ 1),}
 {b_{k} + m_{2,k} + M_{3} (k,s + 1);}
  \hfill \\
 {c_{k} + m_{2,k} + M_{3} (k,s + 1)}; \hfill \\
\end{array}} x_{k}}  \right]. \\
 \end{array}
\end{equation*}

Further, by virtue of the following obvious equalities
$$
N_2(1,s+1) + N_{3} (k,s + 1) = N_{2} (k,s + 1),\,\,1\leq k\leq
s+1, s\in N, $$ $$ m_{2,k} + M_{3} (k,s + 1) = M_{2} (k,s + 1),
1\leq k\leq s+1, s\in N,
$$
we finally find the equality (\ref{e210}). Q.E.D.

\textbf{Remark}. Let $\delta=(\delta_1,...,\delta_n),$ where
$\delta_i\in \{0,1\}$, $1\leq i \leq n$. Let be denote
$$(\alpha):=(\alpha_1,...,\alpha_n),\,\,
(1-\alpha):=(\alpha_1+(1-2\alpha_1)\delta_1,...,\alpha_n+(1-2\alpha_n)\delta_n).$$
Then the following constructive formulas
$$ L^m_{(\alpha)}\left( x_1^{\delta_1(1-2\alpha_1)}\cdot\cdot\cdot x_n^{\delta_n(1-2\alpha_n)}\cdot u \right)= x_1^{\delta_1(1-2\alpha_1)}\cdot\cdot\cdot x_n^{\delta_n(1-2\alpha_n)}L_{{(1-\alpha)}}^m(u)$$
are true. We note, there are $2^n$ equalities.

\section{Fundamental solutions}

Consider equation (\ref{e1}) in $R_{m}^{n +}  .$ Let
$x:=\left(x_1,...,x_m\right)$ be any point and
$x_{0}:=\left(x_{01},...,x_{0m}\right) $ be any fixed point of
$R_{m}^{n +}.$ We search for a solution of (\ref{e1}) as follows:
\begin{equation}
\label{e31} u(x,x_0) = P(r)w(\xi),
\end{equation}
where
$$
\xi = \left( \xi_1,\xi_2,...,\xi_n \right),\,\, \alpha =
\alpha_1+...+\alpha_n - 1 + {\frac{{m}}{{2}}},
$$
$$
P(r) = \left( r^{2} \right)^{ - \alpha}, \,\,\,r^{2} =
{\sum\limits_{i = 1}^{m} {(x_{i} - x_{0i} )^{2}}} ,
$$
$$
r_{k}^{2} = (x_{k} + x_{0k} )^{2} + {\sum\limits_{i = 1,i \ne
k}^{m} {(x_{i} - x_{0i} )^{2}}}, \,\,\xi _{k} = {\frac{{r^{2} -
r_{k}^{2}}} {{r^{2}}}}, \quad k = 1,2,...,n.
$$

We calculate all necessary derivatives and substitute them into
equation (\ref{e1}) :

\begin{equation} \label{e32} {\sum\limits_{k = 1}^{n} {A_{k} \omega _{\xi_{k} \xi
_{k}}} } + {\sum\limits_{k = 1}^{n} {{\sum\limits_{l = k + 1}^{n}
{B_{k,l} \omega _{\xi _{k} \xi _{l}}} } } }    + {\sum\limits_{k =
1}^{n} {C_{k} \omega _{\xi _{k}}} }+ D\omega = 0,
\end{equation}
where
$$A_k=P\sum\limits_{i=1}^m\left(\frac{\partial \xi_k}{\partial x_i}\right)^2, B_{k,l}=2P\sum\limits_{i=1}^m \frac{\partial \xi_k}{\partial x_i}\frac{\partial \xi_l}{\partial x_i}, \,\,k\neq l, k=1,...,n,$$
$$C_k=P\sum\limits_{i=1}^m\frac{\partial^2 \xi_k}{\partial x_i^2}+2\sum\limits_{i=1}^m \frac{\partial P}{\partial x_i}\frac{\partial \xi_k}{\partial x_i}+2P\sum\limits_{j=1}^n \frac{\alpha_j}{ x_j}\frac{\partial \xi_k}{\partial x_j}, $$
$$D=\sum\limits_{i=1}^m\frac{\partial^2 P}{\partial x_i^2}+2P\sum\limits_{j=1}^n \frac{\alpha_j}{ x_j}\frac{\partial P}{\partial x_j}.$$

After several evaluations we find
\begin{equation}
\label{e33}  A_{k} = - {\frac{{4P(r)}}{{r^{2}}}}{\frac{{x_{k}}}
{{x_{0k}}} }\xi _{k} (1 - \xi_{k} ), \end{equation}
\begin{equation}
\label{e34} B_{k,l} = {\frac{{4P(r)}}{{r^{2}}}}\left(
{{\frac{{x_{0k}}} {{x_{k}}} } + {\frac{{x_{0l}}} {{x_{l}}} }}
\right)\xi _{k} \xi _{l} ,\,\,\,k \ne l,\quad l = 1,...,n,
\end{equation}
\begin{equation}
\label{e35}  C_{k} = - {\frac{{4P(r)}}{{r^{2}}}}{\left\{ { -
\xi_{k} {\sum\limits_{j = 1}^{n} {{\frac{{x_{0j}}} {{x_{j}}}
}\alpha _{j}}}   + {\frac{{x_{0k} }}{{x_{k}}} }{\left[ {2\alpha
_{k} - \alpha \xi_{k}}  \right]}} \right\}}, \end{equation}
\begin{equation} \label{e36} D = {\frac{{4\alpha
P(r)}}{{r^{2}}}}{ {  {\sum\limits_{j = 1}^{n} {{\frac{{x_{0j}}}
{{x_{j}}} }\alpha _{j}}}   } } .\end{equation}

Substituting equalities (\ref{e33})-(\ref{e36}) into (\ref{e32})
we obtain the following system of hypergeometric equations of
Lauricella \cite[ p.117]{AP}, which has $2^n$ linearly-independent
solutions \cite[ p.118]{AP}. Considering those solutions, from
(\ref{e31}) we obtain $2^n$ fundamental solutions of equation
(\ref{e1}): \\
\\
$
 1{\left\{ {} \right.} F_A^{(n)} \left[\begin{array}{*{20}c}
a,b_1,...,b_n; \hfill\\ c_1,...,c_n;\hfill\\ \end{array}\xi
 \right],  $\\
$ C_{q}^{1} {\left\{ {{\begin{array}{*{20}c}
 {(x _{1}x_{01})^{1 - c_{1}}  F_A^{(n)} \left[ \begin{array}{*{20}c} a,b_{1} + 1 - c_{1} ,b_{2}
,...,b_{n} ;\hfill\\ 2 - c_{1} ,c_{2} ,...,c_{n} ;\hfill\\ \end{array}\xi  \right]}, \hfill \\
{..............................................................................}
\hfill \\
 {(x _{n}x_{0n})^{1 - c_{n}} F_A^{(n)} \left[\begin{array}{*{20}c} a,b_{1} ,...,b_{n - 1}
,b_{n} + 1 - c_{n} ;\hfill\\ c_{1} ,...,c_{n - 1} ,2 - c_{n} ;
\hfill\\ \end{array}\xi \right]},
\hfill \\
\end{array}}}  \right.}
$ \\
$ C_{n}^{2} {\left\{ {{\begin{array}{*{20}c}
 {(x _{1}x_{01})^{1 - c_{1}}  (x _{2}x_{02})^{1 - c_{2}} F_A^{(n)}
 \left[
\begin{array}{*{20}c} a,b_{1} + 1 - c_{1} ,b_{2} + 1 - c_{2} ,b_{3} ,...,b_{n} ; \hfill\\ 2 -
c_{1} ,2 -
c_{2} ,c_{3} ,...,c_{n} ;\hfill\\ \end{array}\xi  \right]}, \hfill \\
{..............................................................................}
\hfill \\
 {(x _{1}x_{01})^{1 - c_{1}}  (x _{n}x_{0n})^{1 - c_{n}} F_A^{(n)}
 \left[
\begin{array}{*{20}c} a,b_{1} + 1 - c_{1} ,b_{2} ,...,b_{n - 1} ,b_{n} + 1 - c_{n} ;\hfill\\ 2 -
c_{1}
,c_{2} ,...,c_{n - 1} ,2 - c_{n} ;\hfill\\ \end{array}\xi  \right]}, \hfill \\
 {(x _{2}x_{02})^{1 - c_{2}} (x _{3}x_{03})^{1 - c_{3}}
 F_A^{(n)}\left[
\begin{array}{*{20}c} a,b_{1} ,b_{2} + 1 - c_{2} ,b_{3} + 1 - c_{3} ,b_{4} ,...,b_{n}
;\hfill\\ c_{1} ,2 -
c_{2} ,2 - c_{3} ,c_{4} ,...,c_{n} ;\hfill\\ \end{array}\xi \right]},\hfill \\
{..............................................................................}
\hfill \\
 {(x _{n-1}x_{0\,n-1})^{1 - c_{n - 1}} (x _{n}x_{0n})^{1 - c_{n}}  F_A^{(n)}
 \left[
\begin{array}{*{20}c} a,b_{1}  ,...,b_{n - 2} ,b_{n - 1} + 1 - c_{n - 1} ,b_{n} + 1 -
c_{n} ; \hfill\\ c_{1}  ,...,c_{n - 2} ,2 - c_{n - 1} ,2 - c_{n}
;\hfill\\ \end{array}\xi \right]},
\hfill \\
\end{array}}}  \right.}
$ \\$$
.................................................................
$$  $ 1{\left\{{(x _{1}x_{01})^{1 - c_{1}}  \cdot ... \cdot (x _{n}x_{0n})^{1 - c_{n}} F_A^{(n)}\left[ \begin{array}{*{20}c} a,b_{1} + 1 - c_{1} ,...,b_{n} +
1 - c_{n} ;\hfill\\ 2 - c_{1} ,...,2 - c_{n} ;\hfill\\
\end{array}\xi \right]} \right.}, $
\\ where
$$C_n^i=\frac{n!}{i!(n-i)!}, 0\leq i \leq n.$$

It is easy to see that
$$1+C_n^1+C_n^2+...+C_n^{n-1}+1=2^n.$$

The fundamental solutions of equation (\ref{e1}) found above can
be written in a form which is a convenient for further
investigation
\begin{equation}
\label{e37}
 q_k(x,x_0)=\gamma_k r^{-2\alpha}
\prod\limits_{i=1}^n
\left(x_ix_{0i}\right)^{\delta_{ki}(1-2\alpha_i)} \cdot
F_{A}^{(n)}\left[\begin{array}{*{20}c} \alpha+A_k,B_k;\hfill
\\ 2B_k;\hfill\\ \end{array}\xi\right],
\end{equation}
where
$$\delta_k=\left(\delta_{k1},...,\delta_{kn}\right),
\delta_{kj}\in\{0,1\}, \,1\leq j\leq n, \, 1\leq k\leq 2^n,
\,k=1+\sum\limits_{j=1}^n\delta_{kj}\cdot2^{(n-j)\delta_{kj}};$$
$$A_k=\sum\limits_{j=1}^n\left(1-2\alpha_j\right)\delta_{kj},\, B_k=\left(\alpha_1+(1-2\alpha_1)\delta_{k1},...,\alpha_n+(1-2\alpha_n)\delta_{kn} \right),$$
$\gamma_k$ are constants, which will be determined when we solve
boundary-value problems.

\section{Singularity properties of fundamental solutions}

Let us show that the fundamental solutions (\ref{e37}) have a
singularity at $r=0$.

In the case $\delta_1=(0,0,...,0)$ we choose a solution
$q_1(x,x_0)$ and we use the expansion for the hypergeometric
function of Lauricella (\ref{e28}). As a result, a solution
$$ q_1(x,x_0)=\gamma_1r^{-2\alpha}F_{A}^{(n)}\left[\begin{array}{*{20}c}\alpha,\alpha_1,...,\alpha_n;\hfill\\ 2\alpha_1,...,2\alpha_n;\hfill\\ \end{array}\xi\right] $$
can be written as follows
\begin{equation}
\label{e41}
\begin{array}{l}
 q_1(x,x_0)
 = \gamma_1r^{-2\alpha}{\sum\limits_{{\mathop {m_{i,j} = 0}\limits_{(2 \le i \le j \le n)}
}}^{\infty}  {{\frac{{(\alpha)_{N_{2} (n,n)}}} {{{\mathop {m_{2,2}
!m_{2,3} ! \cdot \cdot \cdot m_{i,j} ! \cdot \cdot \cdot m_{n,n}
!}\limits_{(2 \le i
\le j \le n)}}} }}}}  \\
 \cdot{\prod\limits_{k = 1}^{n} {{{{\frac{{(\alpha_{k} )_{M_{2} (k,n)}
}}{{(2\alpha_{k} )_{M_{2} (k,n)}}}
}\left(1-\frac{r_k^2}{r^2}\right)^{M_{2} (k,n)}
F\left[\begin{array}{*{20}c} \alpha + N_{2} (k,n),\alpha_{k} +
M_{2} (k,n);\hfill\\ 2\alpha_{k} + M_{2} (k,n);\hfill\\
\end{array} 1-\frac{r_k^2}{r^2} \right]} }}} .
\\
 \end{array}
\end{equation}

By virtue of formula (\ref{e22}) we rewrite (\ref{e41}) as
$$
 q_1(x,x_0)
 = \frac{\gamma_1}{r^{m-2}}\prod\limits_{k = 1}^{n}r^{-2\alpha_k}_k \cdot f_1\left(r^2,r_1^2,...,r_n^2\right),
$$
where
$$
 f_1\left(r^2,r_1^2,...,r_n^2\right)={\sum\limits_{{\mathop {m_{i,j} = 0}\limits_{(2 \le i \le j \le n)}
}}^{\infty}  {{\frac{{(\alpha)_{N_{2} (n,n)}}} {{{\mathop {m_{2,2}
!m_{2,3} ! \cdot \cdot \cdot m_{i,j} ! \cdot \cdot \cdot m_{n,n}
!}\limits_{(2 \le i \le j \le n)}}} }}}}$$

$$ \cdot{\prod\limits_{k = 1}^{n} {{ {{\frac{{(\alpha_{k}
)_{M_{2} (k,n)} }}{{(2\alpha_{k} )_{M_{2} (k,n)}}}
}\left(\frac{r^2}{r_k^2}-1\right)^{M_{2} (k,n)}
F\left[\begin{array}{*{20}c} 2\alpha_k-\alpha +M_2(k,n)- N_{2}
(k,n),\alpha_{k} + M_{2} (k,n);\hfill\\ 2\alpha_{k} + M_{2}
(k,n);\hfill\\ \end{array} 1-\frac{r^2}{r_k^2} \right]} }}} .
$$

Below we show that $f_1\left(r^2,r_1^2,...,r_n^2\right)$ will be
constant at  $r \to 0$.

For this aim we use an equality (\ref{e21}) and following
inequality
$$N_2(k,n)-M_2(k,n):=\sum\limits_{i=2}^k\left(\sum\limits_{j=i}^n
m_{i,j}-m_{i,k}\right)\geq 0, 1\leq k \leq n \leq m, m>2.$$ Then
we get

\begin{equation}
\label{e42}
 \lim_{r\rightarrow 0}f\left(r^2,r_1^2,...,r_n^2\right)=\frac{1}{\Gamma^n(\alpha)}{\prod\limits_{k = 1}^{n}
 \frac{\Gamma\left(2\alpha_k\right)\Gamma\left(\alpha-\alpha_k\right)}{\Gamma\left(\alpha_k\right)}}.
\end{equation}

 Expressions (\ref{e41})
and (\ref{e42}) give us the possibility to conclude that the
solution $q_1(x,x_0)$ reduces to infinity of the order $r^{2-m}$
at $r \to 0$. Similarly it is possible to be convinced that
solutions $q_k(x;x_0),\,\,k=2,3,...,2^n$ also reduce to infinity
of the order $r^{2-m}$  when $r \to 0$.

It can be directly checked that constructed functions (\ref{e37})
possess properties

\begin{equation}
 {\left. {{\frac{{\partial^{\delta_{kj}} q_{k} \left( {x,x_0}
\right)}}{{\partial x_{j}^{\delta_{kj}} }}}} \right|}_{x_{j} = 0}
= 0, \delta_{kj}\in\{0,1\}, \,1\leq j\leq n, \, 1\leq k\leq 2^n.
\end{equation}

 \end{document}